\documentclass[11pt]{amsart}
\usepackage{amsmath,graphicx}

\theoremstyle{plain}
\newtheorem{theorem}{Theorem}[section]
\newtheorem{lemma}[theorem]{Lemma}

\newtheorem{ther}{Theorem}

\newtheorem{cor}[ther]{Corollary}

\theoremstyle{remark}

\newtheorem*{ack}{Acknowledgement}

\theoremstyle{definition}

\def \R {\mathbf{R}}
\def \Z {\mathbf{Z}}
\def \C {\mathbf{C}}
\def\cp{\mathbf{CP}}
\def\cpbar{\overline{\mathbf{CP}}^{2}}
\def \HH{\mathcal{H}}
\def \WW{\mathcal{W}}
\def\PP{\mathcal{P}}

\def\PPtilde{\tilde\mathcal{P}}
\def\ptilde{\tilde{P}}
\def\cancel#1#2{\ooalign{$\hfil#1\mkern1mu/\hfil$\crcr$#1#2$}}
\def\dirac{\mathrel{\mathpalette\cancel\partial}}
\def\M{\mathcal{M}}

\def\msw{\mathcal{M}^{SW}}
\def\mym{\mathcal{M}^{YM}}

\def\diff{\mathit{Diff}}
\DeclareMathOperator{\SO}{SO}
\DeclareMathOperator{\met}{Met}
\DeclareMathOperator{\spinc}{Spin^c}

\def\SW{\ifmmode{\text{SW}}\else{$\text{SW}$}\fi}
\def\SWW{\ifmmode{\underline{\text{SW}}}\else{$\underline{\text{SW}}$}\fi}
\begin{document}

\baselineskip.525cm
\title{An obstruction to smooth isotopy in dimension 4}
\thanks{The author was partially
supported by NSF Grant 4-50645.}
\author[Daniel Ruberman]{Daniel Ruberman}
\address{Department of Mathematics\newline\indent
Brandeis University \newline\indent
Waltham, MA 02254}
\date\today
\maketitle

\section{Introduction}
It is certainly well-known that a suitable count of solutions to the Yang-Mills or
Seiberg-Witten equations gives rise to invariants of a smooth
$4$-manifold.  In recent years, these invariants have become reasonably
computable, and have led to many advances in $4$-manifold theory.  In this paper,
we use the tools of gauge theory to describe invariants of a diffeomorphism $f$
which depend only on the isotopy class of $f$.  The main application of this study
is the construction of a diffeomorphism
$f$ of a smooth $4$-manifold which is homotopic to the identity, and the proof
that it is not smoothly isotopic to the identity.  It follows from work of
Perron~\cite{perron:isotopy1,perron:isotopy2} and 
Quinn~\cite{quinn:isotopy,freedman-quinn} (see also~\cite{kwasik:isotopy}) that
$f$ is topologically isotopic to the identity.

Formally, the construction of gauge-theoretic invariants of a single
space proceeds as follows.  Under suitable topological hypotheses, the solution
space to the Yang-Mills or Seiberg-Witten equations, denoted $\M(Y;h)$, can be
arranged to be compact and $0$-dimensional.  The parameter $h$
lives in some contractible space $\Pi$ which would be the space of metrics on
$Y$ for Yang-Mills invariants, or the space $\{(g,\delta) \mid *_g\delta= \delta
\}
\subset Met(Y) \times \Omega^2(Y)$ for the Seiberg-Witten invariants. 
One then defines $I(Y)$ to be $\#\M(Y;h)$, where the
symbol `$\#$' denotes an algebraic count of solutions.   While the equations (and
therefore their space of solutions!) involve the parameter
$h$, the parameter space $\Pi$ is connected, so the algebraic count of solutions
is independent of
$h$.  Hence $I(Y)$ depends only on the $C^\infty$ manifold $Y$.

The underlying idea for getting invariants of diffeomorphisms is similar; it
occurs elsewhere in gauge theory~\cite{donaldson:Yang-Mills} and in the context
of symplectic geometry
as well~\cite{bryan-leung:K3,kronheimer:family,seidel:symplectic}.   Suppose
that for some manifold $Y$, the moduli space
$\M(Y;h)$ is formally of dimension $-1$, and hence empty for generic $h \in
\Pi$.  If $f:Y \to Y$ is a diffeomorphism, then we can choose a generic path
$h_t$ from 
$h_0 = h$ to $h_1 =f^*(h)$, and construct the 1-parameter family 
$$
\tilde\M = \bigcup_t \M(Y;h_t)
$$
One then defines an invariant $I(f)$ as the
algebraic count of elements in this $0$-dimensional moduli space.  There are
of course a number of details to check in order to see that $I(f)$ does not
depend on the various choices which were made; this is done below in
section~\ref{defsec}.  The main issue is to see that the invariant is
non-trivial, by computing it in suitable examples.

Our examples can be described rather simply, and reveal the $1$-parameter
invariant to be a sort of suspension of the ordinary YM or SW invariants.  They
arise from the observation that there are many pairs, say $X_0,X_1$ of
non-diffeomorphic $4$-manifolds which are homotopy equivalent, but which
become diffeomorphic after connected sum with $\cp^2$.  Let $Z'$ denote this
connected sum, and note that $Z'$ contains two embedded spheres, say $S_0,S_1$,
of self-intersection $+1$, arising as complex lines in the  two decompositions
$ X_0\#\cp^2 \cong Z' \cong X_1\#\cp^2$.  In typical examples, in fact, $Z'$
decomposes completely into a connected sum of copies of $\cp^2$ and  $\cpbar$.
Assuming this to be the case, a theorem of Wall~\cite{wall:diffeomorphisms} says
that $Z'$ supports enough diffeomorphisms so that we can assume that
$S_0$ and $S_1$ are homologous; of course they are not isotopic.  

Blow up two points on $Z'$ to get $Z = Z' \#_2 \cpbar$, which now contains, in
addition to the $S_i$,  two  $(-1)$-spheres $ E_1,E_2$ which are the the
exceptional curves in $\cpbar$.  Associated to any sphere
$\Sigma$ of self-intersection $-1$ in a $4$-manifold, there is a diffeomorphism
$\rho^\Sigma$ inducing the `reflection in $\Sigma$' on homology:
$$
\rho^\Sigma_*(\alpha) = \alpha + 2(\Sigma\cdot\alpha) \Sigma
$$
From the decomposition of $Z$ as $X_0\# \cp^2 \#_2 \cpbar$, we have
$(-1)$-spheres $S_0+E_1+E_2$ and $S_0-E_1+E_2$, where the `$\pm$' signs refer to
(oriented) connected sums of the spheres in $Z$.  The composition of the
reflection in these spheres is an orientation preserving diffeomorphism of $Z$:  
$$
f_0 = \rho^{S_0+E_1+E_2}\circ\rho^{S_0-E_1+E_2}
$$
Similarly, the decomposition $Z =X_1\# \cp^2 \#_2 \cpbar$ gives rise to a
diffeomorphism 
$$
f_1 = \rho^{S_1+E_1+E_2}\circ\rho^{S_1-E_1+E_2}
$$
The fact that $S_0$ and $S_1$ are homologous implies that $f_0$ and $f_1$ are
homotopic, so that $f = f_1\circ f_0^{-1}$ is homotopic to the identity of $Z$.

The main consequence of the work in this paper is that if there are suitable
degree-$0$  Donaldson invariants $D_{X_j}$ with $D_{X_0}\neq D_{X_1},$ then the
diffeomorphism
$f$ is not smoothly isotopic to the identity.   This follows from the basic
properties of the invariants as outlined above, and the
following calculations.

\begin{ther}\label{maincalc}
The $1$-parameter Yang-Mills moduli space on the manifold $Z$ defines integer
invariants
$D_Z(f_0)$, $D_Z(f_1)$ and $D_Z(f)$, which satisfy:
\begin{enumerate}
\item $D_Z(f_j) = -4D_{X_j}$ for $j=0,1$.
\item $ D_Z(f) = D_Z(f_1) - D_Z(f_0).$
\end{enumerate}
\end{ther}
\begin{cor}
If $D_{X_0}\neq D_{X_1},$ then $f$ is not smoothly isotopic to the identity.
\end{cor}

\begin{ack}
I'd like to thank Tom Mrowka for listening to an earlier, more complicated
version of this work, and for asking the right questions.
\end{ack}
\section{Definition of the invariants}\label{defsec}
The precise nature of the invariants depends on which equations we are using,
so we separate the discussion of the Seiberg-Witten and Yang-Mills equations. 
We will first describe a simple invariant, which corresponds to the degree-$0$
part of the Donaldson polynomial, and then indicate the modifications
necessary to define Seiberg-Witten invariants of a diffeomorphism.  The
Seiberg-Witten invariant, in the version discussed in this paper, is not very
useful, as we shall see.  A more useful Seiberg-Witten invariant, as well as an
extension of the Donaldson invariant to a polynomial invariant, will be
described in subsequent papers~\cite{ruberman:polyisotopy,ruberman:swisotopy}.

\subsection{A simple invariant defined using the Yang-Mills
equations}\label{ymsec} 
We begin by fixing the topological data and some of the notation.  Let $Y$ be an
oriented $4$-manifold, with $b_+^2 > 1$, and suppose that an orientation of
$H^2_+(Y)$ has been chosen.  For simplicity, we will assume that $Y$ is
simply-connected, although this is presumably not necessary.  Fix an $\SO(3)$
bundle $P \to Y$, with non-trivial $w_2$, with the property that the formal
dimension of the ASD moduli space $\mym(P)$ is $-1$.  Note that this implies
that $b^2_+(Y) $ is {\it even}, according to the usual dimension formula. 
Also, for a generic metric $g\in \met(Y)$, the moduli space $\mym(P;g)$ is
empty.  (Here of course we use the non-triviality of $w_2$ and the condition
that $b_+^2(Y)$ be positive to avoid reducibles.  Also, unless there is some
possible confusion with the Seiberg-Witten moduli space, we will drop the $YM$
superscript.)

Although the moduli space is empty for generic $g$, one can consider instead 
a path of metrics $\{g_t \in \met(Y) \mid t \in [0,1]\}$ and the `$1$-parameter
moduli space'
$$
\tilde\M(P;\{g_t\}) = \bigcup_t(\M(P;g_t) \times \{t\})
$$
which evidently has formal dimension $0$.  (The topology on $\tilde\M$ is
as a subspace of $\mathcal{A}/\mathcal{G} \times I$; we will usually omit the
$t$-coordinate.) If $\{ g_t\}$ is a generic path, then $\tilde\M$ will consist
of a finite number of points, which we can count.  Given an orientation on
$H^2_+(X)$, we can assign orientations as well, and get an algebraic count of
points, which we denote by
$ \#\tilde\M(P;\{h_t\})$.
\begin{lemma}\label{ympath}
Suppose that $\M(\ptilde;g_t) = \emptyset$ for
$t= 0,1$.  Then for generic paths $\{ g_t\} \subset \met$ and $\{ g_t'\} \subset
\met$ with $g_0 = g_0'$ and $g_1= g_1'$ we have that 
$$
 \#\tilde\M(P;\{g_t\}) =  \#\tilde\M(P;\{g_t'\})
$$
\end{lemma}
\begin{proof}
Since the space  of metrics is contractible, one can find a $2$-parameter family
$g_{s,t} \in \met$ with $g_{0,t}= g_t$ and $g_{1,t} = g_t'$.   The hypothesis
that $b_+^2 > 2$ implies that for a generic such family, there will be no
reducible solutions to the ASD Yang-Mills equations defined using
the metric ${g_{s,t}}$.  Hence, for a generic family, the $2$-parameter moduli
space $ \Tilde{\Tilde{\mathcal{M}}}$ will be a compact oriented $1$-dimensional
cobordism with boundary $\tilde\M(P;\{g_t'\}) -
\tilde\M(P;\{g_t\})$.
\end{proof}
As a consequence, we have an integer invariant
$D_Y(g_0,g_1)$ defined on pairs of generic metrics on $Y$; note that
interchanging $g_0$ and $g_1$ reverses the orientation of the $1$-parameter
moduli space, and hence changes the sign of $D$.  As
in~\cite{donaldson:Yang-Mills} a similar definition can be given of polynomial
invariants defined using higher-dimensional moduli spaces.  

Suppose now that $f$ is an orientation-preserving diffeomorphism of $Y$, with
the property that $f^*w_2(P) = w_2(P)$.  Since $f$ has degree one, it
automatically pulls back $p_1(P)$ to itself, and hence $f^*(P) \cong P$.  It
is worth remarking that while there is no canonical isomorphism between
$f^*(P)$ and $P$, any two isomorphisms differ by composition with an element of
the gauge group of $P$.   Let $\alpha(f)$ be the spinor norm of $f$, i.e., $\alpha(f) =
\pm 1$ depending on whether $f^*$ preserves or reverses orientation on
$H^2_+$.    Recall that an orientation of the moduli space is determined by the
homology orientation, and the choice of an integral lift $c \in H^2(Y;\Z)$ for
$w_2(P)$.  Let $\beta(f) = (-1)^{(\frac{f^*c-c}{2})^2}$, then it is readily
checked that $\beta(f)$ does not depend on the choice of $c$.   We will say
that $f$ preserves or reverses orientation on the moduli space depending on
whether the quantity $\alpha(f)\beta(f) = \pm 1$.  (This terminology is explained
by Corollary 3.28 of~\cite{donaldson:orientation}.) 

Let $g_0$ be a metric for which the moduli space is empty; we propose to define
an invariant $D_Y(f) = D_Y(g_0,f^*g_0)$, i.e., the count of points in a
$1$-parameter moduli space $\tilde\M(P;g_t)$ where $g_1=f^*g_0$.  
Lemma~\ref{ympath} says that this is independent of the chosen path; a
more subtle point is that  the count is also independent of the initial metric
$g_0$, so long as the final metric $g_1$ is $f^*g_0$.
\begin{theorem}\label{ymwelldef} Suppose that $b_+^2(Y) > 2$.Let $g_0$ and
$k_0$ be generic metrics on $Y$, and let $g_t,k_t \in \met(Y)$ be generic paths
with $g_1=f^*g_0$ and $k_1 = f^*k_0$.  If $\alpha(f)\beta(f) = 1$, then the
algebraic number of points in the respective $1$-parameter moduli spaces are
the same:
$$
\#\tilde\M(P;\{g_t\}) = \#\tilde\M(P;\{k_t\}) 
$$
If  $\alpha(f)\beta(f) = -1$, then the same statement holds modulo $2$.
\end{theorem}
Because of the theorem, we can unambiguously define
$$
D_Y(f) = \#\tilde\M(P;\{g_t\}) \in \Z\ \text{or}\ \Z_2
$$
(depending on the sign of $\alpha(f)\beta(f)$) using any generic $g_0 \in \met(Y)$.
\begin{proof}[Proof of Theorem~\ref{ymwelldef}]
The key observation is that in general, for $g \in \met$, there is an
isomorphism of moduli spaces 
$$
f^* : \M(P,g) \xrightarrow{\cong} \M(P,f^*g)
$$
This means that $f^*$ is a homeomorphism, which is an isomorphism on the level
of deformations, or more formally of locally ringed spaces.  (As remarked above,
$f^*$ is only well-defined after we have divided by gauge-equivalence, because
it involves a choice of isomorphism between $f^*P$ and $P$.) 

Start by choosing a generic path $K_{0,t}$ from $g_0 = K_{0,0}$
to $k_0 = K_{0,1}$, and note that $K_{1,t} = f^*K_{0,t}$ is thus necessarily
generic.  Now take arbitrary paths $g_s$ and $k_s$, and define $K_{s,0} = g_s$
and $K_{s,1} = k_s$.  Thus we have a well-defined loop in the space of metrics,
and as in the proof of Lemma~\ref{ympath}, can fill it in with a generic
$2$-parameter family $K_{s,t}$.  

Now consider the $2$-parameter moduli space
$$
\Tilde{\Tilde{\mathcal{M}}}  = \bigcup_{s,t} \M(P,K_{s,t})
$$
which again is a compact $1$-manifold with boundary
$\partial\Tilde{\Tilde{\mathcal{M}}}$.  (For the compactness, we need that
there are no reducibles in a 2-parameter family, which is why we require that
$b_+^2 > 2$.) By construction,
$f^*$ induces an orientation preserving diffeomorphism between the part of the
boundary lying on the path
$K_{0,t}$ to the part lying on the path $K_{1,t}$. Put another way, the the
algebraic number of points in the $1$-parameter moduli space
$\bigcup_t\M(P,K_{0,t})$ must be the same as the number of points in 
$\bigcup_t\M(P,K_{1,t})$.   Hence the same must be true for the algebraic
count of the points in $\bigcup_s \M(P,K_{s,0})$  and $\bigcup_s
\M(P,K_{s,1})$.  This proves the theorem; see figure 1 for an illustration.
\end{proof}
\begin{center}
\includegraphics{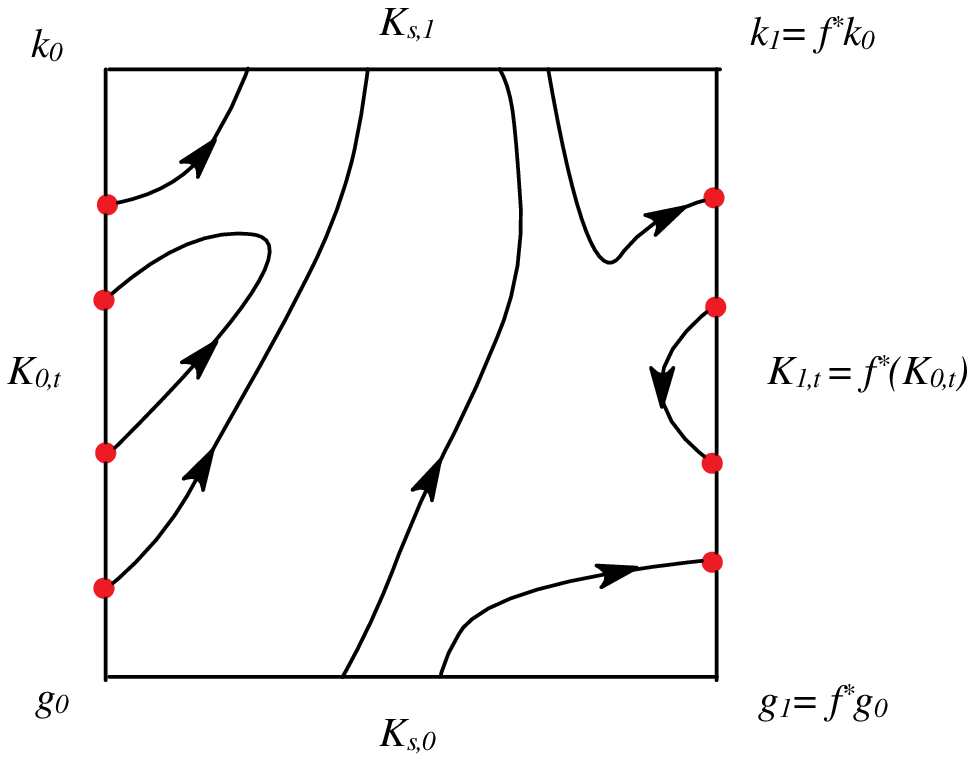}\\
Figure 1
\end{center}

\subsection{Invariants defined via the Seiberg-Witten equations}\label{swsec}

The definition of an invariant of diffeomorphisms based on the Seiberg-Witten
equations is formally similar to that given above for the Yang-Mills equations. 
The main difference is that while the condition that $f$ preserves the isomorphism
class of an $\SO(3)$ bundle is fairly weak, the condition that $f$ preserves a
$\spinc$ structure is much more restrictive.  In this section, we discuss how to
carry through the definition; those proofs which are essentially parallel to ones
in section~\ref{ymsec} are omitted or only sketched.  Seiberg-Witten invariants
of $f$ can be defined without the restriction that $f$ preserve a $\spinc$
structure; this will be treated in more detail in~\cite{ruberman:swisotopy}.

Suppose that $(Y,g)$ is a Riemannian $4$-manifold, with principal frame bundle
$P = P(g) \rightarrow Y$.  A $\spinc$ structure on $Y$, from the point of view
of principal bundles, is an $S^1$ bundle $\ptilde \rightarrow P $ which
restricts on each fiber to the bundle $\spinc(4) \rightarrow \SO(4)$.  Although
the $\spinc$ structure is described in terms of the metric on $Y$, it is
possible to show that a $\spinc$ structure for one metric induces one for all
metrics, as follows.  Note that there is a `universal' $\SO(4)$ bundle $\PP \to
Y \times \met(Y)$ whose restriction to any slice $Y \times \{g\}$ is the frame
bundle of $g$.  Since the space $\met(Y)$ is contractible, the bundle $\ptilde
\to P(g)$ induces in a canonical way a bundle $\PPtilde \to Y \times \met$ which
gives a
$\spinc$ structure for each metric.  Under the projection $\Pi \to \met(Y)$, the
universal bundle pulls back to a bundle which we continue to denote $\PPtilde
\to Y \times \Pi$.  The spin representations give bundles $\mathcal{W}^\pm$
over $Y \times \Pi$ whose restrictions to $Y \times \{h\}$ are the ordinary
spin bundles $W^\pm(g)$.

Let $f$ be an orientation preserving diffeomorphism of $Y$, and note that if
$\delta $ is a $g$ self-dual form, then $f^*\delta$ is $f^*g$ self-dual.  Hence
$f$ has a natural action on $\Pi$:
$$ h = (g,\delta) \rightarrow f^*h = (f^*g,f^*\delta).
$$
Associated to any $\spinc$ structure and $h\in \Pi$, we have the (perturbed)
Seiberg-Witten equations for a connection $A$ on $\det(W^+)$ and spinor
$\phi \in \Gamma(W^+(g))$
\begin{equation*}\tag{$SW_{h}$}\label{swh}
\begin{split}
F_A &= q(\phi) + i\delta\\
\dirac_A(\phi) &= 0
\end{split}
\end{equation*}
Denote the solution space, modulo gauge equivalence, by $\M(\ptilde;h)$ (or
perhaps $\msw(\ptilde;h)$ if it's not clear what equations are being
discussed.) Let us assume that $b^2_+(Y)$ is even, and that for the $\spinc$
structure $\ptilde$, the formal dimension of $\msw(\ptilde)$ is
$-1$.  Thus for generic $h \in \Pi$, the moduli space is empty.  Once again, we
consider the $1$-parameter moduli space
$$
\tilde\M(\ptilde;\{h_t\}) = \bigcup_{t} \M(\ptilde(h_t);h_t)
$$
which evidently has formal dimension $0$.  Given an orientation on $H^2_+(Y)$,
we can assign orientations to the points, and get an algebraic count of points,
which we define to be $\SW(\ptilde;\{h_t\})$.  We then have the analogue of
Lemma~\ref{ympath}, leading to an invariant of pairs of generic elements $h_0,h_1
\in \Pi$.
\begin{lemma}\label{swpath}
Suppose that $\M(\ptilde;h_t) = \emptyset$ for
$t= 0,1$.  Then for generic paths $\{ h_t\} \subset \Pi$ and $\{ h_t'\} \subset
\Pi$ with $h_0 = h_0'$ and $h_1= h_1'$ we have that 
$$
 \SW(\ptilde;\{h_t\}) =  \SW(\ptilde;\{h_t'\})
$$
\end{lemma}

Now suppose that $f:Y \xrightarrow{\cong} Y$ is a diffeomorphism which
preserves the orientation of $Y$ and acts on the orientation of $H^2_+(Y)$ with
sign $\alpha(f)$, and that $\ptilde$ is a $\spinc$ structure such that
$f^*(\PPtilde) \cong \PPtilde$.  Suppose that $h_0\in \Pi$ is generic, so that
$\M(\ptilde;h_0) $ is empty.   Let $h_1 = f^*(h_0)$, and take a generic path
$\{h_t\}$ from $h_0$ to $h_1$.   By analogy with Theorem~\ref{ymwelldef}, and its
proof, we have:
\begin{theorem}\label{swwelldef}  If $k_0 \in \Pi$ is also a generic point, and
$\{ k_t\} $ is a generic path from $k_0$ to $k_1 = f^*k_0$, then 
$$
 \SW(\ptilde;\{h_t\}) =  \SW(\ptilde;\{k_t\})
$$ 
as integers if $\alpha(f) = +1$ and modulo $2$ if $\alpha(f) = -1$.
\end{theorem}
Because of the theorem, we can unambiguously define 
$$ \SW(f,\ptilde) = \SW(\ptilde;\{h_t\}) \in \Z\ \text{or}\ \Z_2 $$
using any generic $h_0 \in \Pi$.

The condition that $f$ preserve the $\spinc$ structure can be relaxed somewhat. 
In fact, this condition was used only at one point in the proof, namely in
identifying 
$\bigcup_t\M(\ptilde,K_{0,t})$ with
$\bigcup_t\M(\ptilde,K_{1,t})$.  Recall~\cite{morgan:swbook} that there is a
natural involution
$J$ on the space of $\spinc$ structures, essentially given by replacing the spin
bundle by its complex conjugate.   The principal $\spinc_4$ bundle corresponding
to a
$\spinc$ structure $\ptilde$ will be denoted $-\ptilde$.  A fundamental
observation is that $J$ behaves in a natural way with regard to the
Seiberg-Witten equations, and in particular gives an isomorphism $J$ of moduli
spaces
$$
\msw(\ptilde;h) \cong \msw(-\ptilde;-h)
$$
Here, if $h = (g,\delta)$, then $-h$ will be the pair $(g,-\delta)$.  The
orientations on $\msw(\ptilde;h) $ and $ \msw(-\ptilde;-h)$
differ~\cite[\S6.8]{morgan:swbook} by a sign of
$(-1)^{\epsilon}$, where $\epsilon = 1 +b_1(Y) + b^2_+(Y) + \text{ind}_{\bf
C}(\dirac)$.   The assumption that $H_1(Y) = 0$ and that the moduli space is
$-1$ dimensional shows that in the case at hand, $\epsilon =
\frac{b^2_+(Y)}{2} + 1$.    Let
$\alpha(f)$ denote the spinor norm of
$f$, i.e.~
$\alpha(f) =
\pm 1$ depending on whether $f$ preserves or reverses the homology orientation.

Suppose that $f$ is an orientation preserving diffeomorphism, such that
$f^*\ptilde \cong -\ptilde$.  Then we have the following variation of
Theorem~\ref{swwelldef}.
\begin{theorem}\label{swbar}
Let $f$ be an orientation preserving diffeomorphism, such that $f^*\ptilde \cong
-\ptilde$.  Let $h_0,k_0$ be generic points in  $\Pi$, and suppose that
$\{h_t\}$ is a generic path from $h_0$ to $h_1 = f^*h_0$, with similarly chosen
generic path $\{k_t\}$.  Then 
$$
 \SW(\ptilde;\{h_t\}) = \SW(\ptilde;\{k_t\})
$$ 
where the equality is over the integers if  $(-1)^{\epsilon}\alpha(f) = 1$ and
is only modulo $2$ if $(-1)^{\epsilon}\alpha(f) = -1$.
\end{theorem}
The proof is essentially the same as that of Theorem~\ref{swwelldef}, using the
composition of diffeomorphisms
$$
 \M(\ptilde;h) \xrightarrow{f^*} \M(-\ptilde,f^*h) \xrightarrow{J}
\M(\ptilde,-f^*h) 
$$
To see that the assumption that  $(-1)^{\epsilon}\alpha(f) = -1$
means that count of solutions is well-defined only modulo $2$, the reader
should imagine reversing all the arrows on the curves in Figure~1 which meet
the left-hand edge (where
$t=0$.)

\subsection{Some simple properties of the invariants}
The invariants we have defined have many similarities with the ordinary
Donaldson or Seiberg-Witten invariants; for instance they change sign when the
homology orientation is reversed.  In addition, they satisfy gluing properties
analogous to those which hold for the ordinary invariants, the simplest of
which we will exploit in section~\ref{wallsec} to compute a non-trivial
example.  We will not develop this systematically, but instead record here some
basic facts which will be needed in the computation.  To avoid repeating
essentially identical argument for the different sorts of invariants defined in
the preceding sections, we will denote any of $D_Y(f)$,
$D_Y(f,\Sigma_1,\ldots,\Sigma_d)$ or $\SW(f,\ptilde)$ simply by $I(f)$. 
Similarly,  $\Pi$ will denote either the space of metrics on $Y$ (if we are
discussing Yang-Mills invariants) or the space of metrics and self-dual forms
(for the case of Seiberg-Witten invariants), and $\M(h)$ will denote a moduli
space.   The first simple fact is that 
$I(f)$ is additive under composition.  To put it in another way, the set of
diffeomorphisms for which a given invariant is defined is a subgroup
of $\diff^+(Y)$, and $I$ gives a homomorphism to  $\Z$ or $\Z_2$. 
\begin{lemma}\label{additive}
Suppose that $f$ and $f'$ are orientation preserving diffeomorphisms of $Y$, and
that the invariants $I(f)$ and $I(f')$ are defined.  Then
\begin{enumerate}
\item $I(f'\circ f)$ is defined, and equals $I(f') + I(f)$.  
\item  $I(f^{-1})$ is defined and equals  $-I(f)$.
\end{enumerate}
\end{lemma}
\begin{proof}
We make use of the fact that initial points
$h_0\in \Pi$ and paths $\{h_t\} \in \Pi$ may be chosen arbitrarily.  Let $h_0$ be 
a generic point in $\Pi$, and choose a generic path $h_t$, $t\in [0,1]$ with $h_1
= f^*h_0$, so that $\# \M(\{h_t\}) = I(f)$.  Choose a generic path, say $h_t'$, $t
\in [0,1]$, with $h_0' = h_1$, so that $\# \M(\{h_t'\}) = I(f')$.  Putting the
two together gives a path $h*h'$, which calculates $I(f'\circ f)$.  But
manifestly the (algebraic) number of points in the corresponding
$1$-parameter moduli space is
$\# \M(\{h_t\}) + \# \M(\{h_t'\}) = I(f') +I(f)$.

The second part follows by considering the initial point in $\Pi$ for $f^{-1}$
to be
$f^*h_0$, and using the path $h_{1-t}$.  The reversal in the path changes the
orientation of the $1$-parameter moduli space, which accounts for the minus sign.
\end{proof}

The other important principle to note is that
$I(f)$ depends only on the isotopy class of $f$.
\begin{lemma}\label{isotopy}
Suppose that $f_0$ and $f_1$ are smoothly isotopic diffeomorphisms of a
$4$-manifold $Y$, and that one of the invariants $I(f)$ is defined. Then
$I(f_0)=I(f_1)$.
\end{lemma}
\begin{proof}
Note that if $f_0$ is isotopic to $f_1$, then $f_1\circ f_0^{-1}$ is isotopic to
the identity.  So by the previous lemma, and the obvious remark that the
identity has trivial invariant, it suffices to show that $I(f_1) = 0$ if $f_1$ is
isotopic to the identity.  Let $f_t$ be an isotopy from the identity to
$f_1$.  Choose 
$h_0\in \Pi$  for which the moduli space $\M(Y;h_0)$ is empty, and consider the
path $h_t = f_t^*h_0$.  As remarked in the proof of Theorem~\ref{ymwelldef},
for each $t\in [0,1]$ the moduli space
$\M(Y;h_t)$ is isomorphic to  $\M(P;h_0)$, i.e., is empty.  Hence $I(f_1) =
0$.
\end{proof}

\section{Smoothly non-isotopic diffeomorphisms}
In this section, we give the detailed construction of the diffeomorphism
described in the introduction, and use the invariant $D(f)$ to show that
it is not smoothly isotopic to the identity.    As described in the introduction,
the diffeomorphism is a composition of four `reflections' on a manifold $Z$,
where each reflection is associated to a certain $2$-sphere embedded in $Z$.  So
as an initial step, we discuss the $1$-parameter gauge theory associated to a
single reflection.

The manifold $Z$ will be written as $X_0\# \cp^2\#_2\cpbar$, with some
topological hypotheses on $X_0$ needed for the calculation to succeed.  These
depend on whether one wants to calculate $D(f)$ or $\SW(f)$,
and somewhat surprisingly are considerably more restrictive in the
Seiberg-Witten case.  In either case, we assume that $H_1(X_0) = 0$, and that 
$b^2_+$ is odd and greater than $1$.  In addition, for the Yang-Mills version, we
require
\begin{equation*}\tag{$*_{YM}$}\label{ymhyp}
\begin{split}
{}&\text{There is an $\SO(3)$ bundle $P_0\to X_0$ with}\\
{}&\text{$w_2(P_0) \neq 0$ and
$\dim(\mym(P_0)) = 0$.}
\end{split}
\end{equation*}
The requirement for the Seiberg-Witten version is discussed in
section~\ref{swversion}.

Let $N = \cp^2 \#_2\cpbar$, and let $S, E_1,E_2  $ be the obvious embedded 
$2$-spheres of self-intersection $\pm 1$, and let $\Sigma_\pm = S\pm E_1+E_2$ be
spheres of square $-1$.  For $\Sigma = \Sigma_+$ or $\Sigma_-$, there is a
diffeomorphism $\rho^\Sigma$ of
$N$, inducing in homology a reflection in $\Sigma$, in other words,
$\rho^\Sigma_*(x) = x + 2(x\cdot \Sigma)\Sigma$, where `$\cdot$' is the
intersection product.  (In cohomology,
$(\rho^\Sigma)^*$ is given by the same formula, where  $\Sigma$ is replaced by
$\sigma =  PD(\Sigma)$ and
$\cdot$ is replaced by
$\cup$.)  In particular, letting $s,e_1,e_2$ be the Poincar\'e duals of $S,
E_1,E_2$, we have $(\rho^\Sigma)^*(s) = s + 2(s\cup \sigma)\sigma= 3s
+2(e_2 \pm e_1)$.  Hence $\rho^\Sigma $ is orientation preserving on $H_+^2(N)$. 
Choose the orientation in which $s$ is in the positive cone.

Fix a point $p$ in $N$, and choose specific
diffeomorphisms $\rho^\Sigma$ which are equal to the identity in some fixed
neighborhood of $p$.  For any $4$-manifold $X_0$, such a diffeomorphism glues
together with the identity map of $X_0$ to give a diffeomorphism $f^\Sigma$ on the
connected sum $Z = X_0\# N$.  Let $L (=L(1,1,1))$ be the complex line bundle
over $N$ with first Chern class $s+ e_1 + e_2$. The bundle
$P_0$ may be glued together with
$L \oplus \R$ to give an $\SO(3)$ bundle $P_Z \to Z$, which will be the basis
for our computation.  

Any single diffeomorphism $f^\Sigma$ reverses the orientation of the moduli
space, and so has only $\Z_2$ invariants.  In fact, it turns out (see the proof
of Theorem~\ref{singlecalc} and the discussion after Theorem~\ref{swcalc}) that
all of these
$\Z_2$ invariants are trivial.  To get a more interesting computation, we
compose two reflections, to get $f_0 = f^{\Sigma_+}\circ f^{\Sigma_-}$.  At the
moment we do not know how to define an interesting Seiberg-Witten invariant of
$f_0$, so we will concentrate mainly on the calculation of a Donaldson invariant
$D_Z(f_0)$. In the next section, we will verify the first part of
Theorem~\ref{maincalc} of the introduction.
\begin{theorem}\label{singlecalc}
Suppose that $P_0 \to X_0$ is an $\SO(3)$ bundle with non-trivial $w_2$, so that
the $0$-dimensional Donaldson invariant $D_{X_0}$ is defined.  Then the
invariant
$D_Z(f_0) \in \Z$ is defined, and $D_Z(f_0) = -4D_{X_0}$.
\end{theorem}

It is worth noting that this theorem depends only on the homotopy type of $N$,
and on the maps induced in cohomology by the restriction of $f_0$ to $N$.

\subsection{Wall-crossing and the $1$-parameter invariant}\label{wallsec}

The calculation of the invariant $D_Z(f_0)$ depends on a simple wall-crossing
argument in the style of~\cite{donaldson:dolgachev,kotschick:so3}.  We give the
details of this argument in this subsection.  First we need to review the chamber
structure associated to reducible connections on a particular $\SO(3)$ bundle
over $N$. 

On $N$, consider the $\SO(3)$ bundle $P$ with associated $\R^3$ bundle 
$L\oplus\R$, where $L$ is the complex line bundle with $c_1 = s+e_1+e_2$. 
Reductions of $P$ to an $\SO(2)$ bundle are in one-to-one correspondence with
cohomology classes $(a,b,c) \in H^2(N;\Z)$ where $c$ is positive, $a,b,$ and $c$
are odd and satisfy $c^2-a^2-b^2= -1$.   For any metric $g$ on $N$, there
is a unique harmonic form $\omega_g \in \Omega_+^2(N)$ with 
\begin{equation}\label{hyper}
\int_N \omega_g \wedge \omega_g = 1 \quad
\text{and}\quad 
[\omega_g]\cup s > 0
\end{equation}
Denote by $\HH\subset H^2(N;\R) $ the set of cohomology classes
satisfying~\eqref{hyper}. The line bundle corresponding to a triple $(a,b,c)$
admits a
$g$-ASD connection if and only if $\omega_g$, as an element of $H^2(N;\R)$, lies
on the `wall'
$$W(a,b,c) = \{(x,y,z)\in \HH \mid cz -ax-by = 0\}.$$
It is readily verified that distinct walls are in fact disjoint
(cf.~\cite[Chapter II]{friedman-morgan:complex-surfaces-I}; a picture of the
 walls may be found below in Figure 2.)  The
walls are transversally oriented, using the following convention.    Let
$\epsilon(a,b,c)
$ be the sign with which the orientation of $\M(P_N)$ induced by the reduction
to $L(a,b,c)$ compares to the given one, i.e. that induced by the reduction to
$L(1,1,1)$.  According to~\cite{donaldson:orientation}, this is given by 
$$
\epsilon(a,b,c)= (-1)^{[(\frac{c-1}{2})^2 - (\frac{a-1}{2})^2 
-(\frac{b-1}{2})^2]}
$$
The transverse orientation of $W(a,b,c)$ is chosen so that the positive
side consists of those $(x,y,z)$ with
$\epsilon(a,b,c) (cz-ax-by) > 0$.  Let us write $\WW $ for the union of the
$W(a,b,c)$; it is a transversally oriented, locally finite, proper submanifold
of $\HH$.

Write  $Z = X\# N$; eventually $X$ will be either $X_0$ or $X_1$.  We
assume that we have an $\SO(3)$-bundle $P_X\to X$ with non-trivial $w_2$,
for which the moduli space is formally $0$-dimensional.  Consider a metric
$g^X$ which is Euclidean in a neighborhood of the point $p$ where the connected
sum is made.  The complement of $p$ has a neighborhood which is conformally
the  cylinder $S^3 \times [-T,1)$.  By a small perturbation away from this tube,
it may be assumed that the moduli space $\mym(P_X;g^X)$ is smooth; the fact that
$w_2 \neq 0$ readily implies that it is also compact.  Choosing an orientation
of $H^2_+(X)$, together with an integral lift $c_X$ of $w_2(P_X)$ yields the
Donaldson invariant $D_{X}$.   

The same data gives a homology orientation for $Z$; it is the
homology orientation $X$, followed by  the choice of generator $s\in
H^2(\cp^2) \cong H^2_+(N)$.  We fix the integral lift $c_X \oplus (s+e_1+e_2)
\in H^2(X) \oplus H^2(N)$ for $w_2(P_Z)$.  Note that 
$$
f^*(c_X \oplus (s+e_1+e_2)) =  c_X \oplus (-3s+ e_1-3e_2)
$$
so that the data needed to define the integral invariant $D_Z(f_0)$ is
present.

Choose a metric $g_0^N$ on $N$ with a similar structure near the connected sum
point, and glue together to get a metric $g_0$ on $Z$, so that $Z$ contains a
long tube (of length $2T$).  By a perturbation of $g_0^N$, supported away from
$p$, we can assume that the harmonic form $\omega_{g_0^N}$ (or $\omega_0$ for
short) does not lie on any wall.  Furthermore, since the dimension of $\mym(P_N)$
is $-2p_1(L\oplus \R) -3(b^2_+(N) +1) = -4$, we can assume that 
$$
\mym(P_N;g_0^N) = \emptyset
$$
Suppose that $g_1^N$ is a second metric which satisfies the same genericity
properties and has the same behavior near $p$; of course $g_1^N$ defines a
self-dual harmonic form $\omega_1$.  This metric may also be glued to the metric
$g^X$ to obtain a metric
$g_1$ on
$Z$.  We want to evaluate the
$1$-parameter Donaldson invariant $D_Z(g_0,g_1)$, at least in the limit when
the tube length $T$ is large.  To this end, choose any path in $\gamma\subset \HH$
from $\omega_0$ to
$\omega_1$ which is transverse to all of the walls.
\begin{lemma}\label{wallcross}
If $T$ is sufficiently large, then
$D_Z(g_0,g_1) = 2(\gamma \cdot \WW) D_X$.
\end{lemma}
\begin{proof}
We begin by remarking that the intersection number in the statement of the lemma
is indeed well-defined, i.e., independent of the path in $\HH$ provided the
endpoints are fixed.  Thus, we are free to choose any path of metrics on $Z$;
our path will be of the form $g^X \# g^N_t$, where the metrics on $N$ are
chosen as follows.  Make an initial choice of path so that the associated path
of harmonic forms is transverse to $\WW$.  Consider a compact subinterval $J$ of
$[0,1]$ for which the associated harmonic forms do not lie on $\WW$.   Because
the parameterized moduli space $\mym(P_N,g^N_t)$ is formally of dimension $-3$,
a small perturbation (relative to the endpoints of
$J$) will ensure that $\mym(N;J) = \emptyset$.  This fails to hold for intervals
whose associated path of harmonic forms meet $\WW$, because one cannot perturb
away the reducible connections. 

Each intersection point of $\gamma$ and $\WW$ thus contributes to
$D_Z(g_0,g_1)$, and we need to show that the contribution is given by $\pm 2
D_X$ when the intersection is transverse, and to evaluate the sign. 
Consider, for $t \in (-\epsilon,\epsilon)$, a path of metrics so that the
associated path of harmonic forms crosses the wall $W(a,b,c)$ transversally 
at $t=0$.  For each $t$, there is a unique unitary connection $A^N_t$
on $L$ such that $\frac{1}{2\pi i}F_{A_t}$ is $g^N_t$-harmonic.  The assumption
that the path $\omega_{G^N_t}$ crosses the wall transversally at $t=0$ means
precisely that $A^N_0$ is anti-self-dual, and that the function 
\begin{equation}\label{sdproj}
h(t) = (a e_1 + b e_2 + c s) \cup \omega_{g^N_t} = \frac{1}{2\pi i}F_{A_t}  \cup
\omega_{g^N_t} = \frac{1}{2\pi i}F_{A_t}^+  \cup
\omega_{g^N_t}
\end{equation}
has non-zero derivative at $0$. 

One can glue the whole path $A^N_t$ to any given ASD connection $A^X$ on $P_X$,
and consider the problem of deforming the resulting path to ASD connections on
$P_Z$. The only possible solutions (for sufficiently large tube length) are
when $t$ is near zero, and may be described in terms of the Kuranishi picture
for the $1$-parameter gluing problem. 
The Kuranishi picture (compare~\cite{donaldson-kronheimer,kronheimer-mrowka:thom})
for the gluing problem is an
$\SO(2) =\text{Stab}(A^N_0)$-equivariant map 
$$
\psi : (-\epsilon,\epsilon) \times \HH^1(A^N_0)   \times  \HH^1(A^X) 
\times \SO(3) \rightarrow
\HH^2_+(A^N_0) \times \HH^2_+(A^X)
$$
so that the part of the moduli space $\M(P_Z;\{g_t\})$ for $t \in
(-\epsilon,\epsilon)$ is modeled by $\psi^{-1}(0)/\SO(2)$.
 
Since $A^X$ is a smooth point in the $0$-dimensional space $\mym(P_X)$, it
has vanishing cohomology groups $\HH^1(A^X)$ and $\HH^2_+(A^X)$.  On the $N$ side,
the cohomology groups split as $\HH^*(L) \oplus \HH^*(N;\R)$, where the first summand
gets a complex structure from $L$.    The index calculation shows that
$\HH^1(A^N) \cong \HH^1(L)
\cong \C^r$ and $\HH^2_+(A^N) \cong \HH^2_+(L) \oplus \HH^2_+(N;\R )\cong \C^{r+1} 
\oplus \R$ for some $r \ge 0$.  In fact~\cite{freed-uhlenbeck} by a
small perturbation of the whole path of metrics on $N$, one can arrange that
$r=0$, which we do for convenience.  Hence 
$\psi$ is an $\SO(2)$ equivariant map
$$
(-\epsilon,\epsilon)  \times \SO(3) \rightarrow \C \times \R
$$
The projection
of $\psi$ onto the last component (ie the one corresponding to $\HH^2_+(N;\R)$) is
given, at the point $(0,0)$ by the map $h$ described in equation~\eqref{sdproj},
up to first order.  Thus the local contribution of $A^X$ is given by counting the
zeros of an $\SO(2)$ equivariant map
$$
\psi_0: \SO(3) \rightarrow \C 
$$ 
Such a map is equivalent to a section of the complex line bundle
$$
\SO(3) \times_{\SO(2)} \C \rightarrow \SO(3)/\SO(2) \cong S^2
$$
which is readily seen to be the tangent bundle of $S^2$.  Since the Euler class
of the tangent bundle is $2$, each point $A^X \in \M(P_X)$ contributes $\pm 2$
to the $1$-parameter invariant.  To compute the sign, note first that the
determinant bundle for $d^* + d^+$ has an orientation coming from the reduction
into $L(a,b,c) \oplus \R$.  As remarked before, this differs from the
orientation arising from the original splitting into $L(1,1,1) \oplus \R$, with
sign given by $\epsilon(a,b,c)$.  Hence each intersection point contributes $2
\epsilon$ times the sign of $h'(0)$, which is exactly the intersection number
with $W(a,b,c)$, using the prescribed transverse orientation.
\end{proof}

With this basic step in hand, we can compute the invariant for the
diffeomorphism $f_0$ described above.
\begin{proof}[Proof of Theorem~\ref{singlecalc}]
Observe that each of the diffeomorphisms $f^{\Sigma_\pm}$ is orientation
preserving, so that to see that there is a well-defined invariant $D_Z(f_0)$,
it suffices to check the effect of $f_0$ on $c_1(L) = s + e_1+ e_2$.  It is
readily verified that $f^*_0(s + e_1+ e_2) = 5s + e_1+ 5e_2$, which shows that
$f_0$ preserves $P_0$, and also that the action on the moduli space is
orientation preserving.  Hence the integer invariant $D_Z(f_0)$ is defined. 
According to the calculation in Lemma~\ref{wallcross}, we only have to check,
for any
$\omega_0 \in \HH$, the intersection number of a path in $\HH$ from $\omega_0$
to
$f^*_0(\omega_0)$ with the walls $\WW$.  But this is easily seen to be be $-2$,
by taking a point such as $(0,0,1)$; see the following figure and discussion.
\end{proof}

A few of the walls are illustrated below, in the upper half-space model of
$\HH$.  The point $(0,0)$ (which corresponds to $(0,0,1)$ in the hyperboloid
$z^2-x^2-y^2 = 1$) is indicated with a small $+$ sign, together with its image
under
$f_0^*$.  The path from $(0,0)$ to its image crosses the walls corresponding to
the reductions $L(1,1,1)$ and $L(1,3,3)$, each with a negative orientation. 
As a check of the basic set-up, a second point $(-1/2,-1/2)$, as is its
image, is indicated by a small circle $\circ$.  The path from
$(-1/2,-1/2)$ to $f_0^*(-1/2,-1/2)$ crosses the walls $W(-1,-1,1)$, $W(1,1,1)$,
$W(1,3,3)$, and $W(7,11,13)$.  The first crossing counts for $+2$, while the
others count $-2$, for a net of $-4$.

\begin{center}
\includegraphics[trim= 1.4in 2.5in 1.4in 2.5in, width=5in]{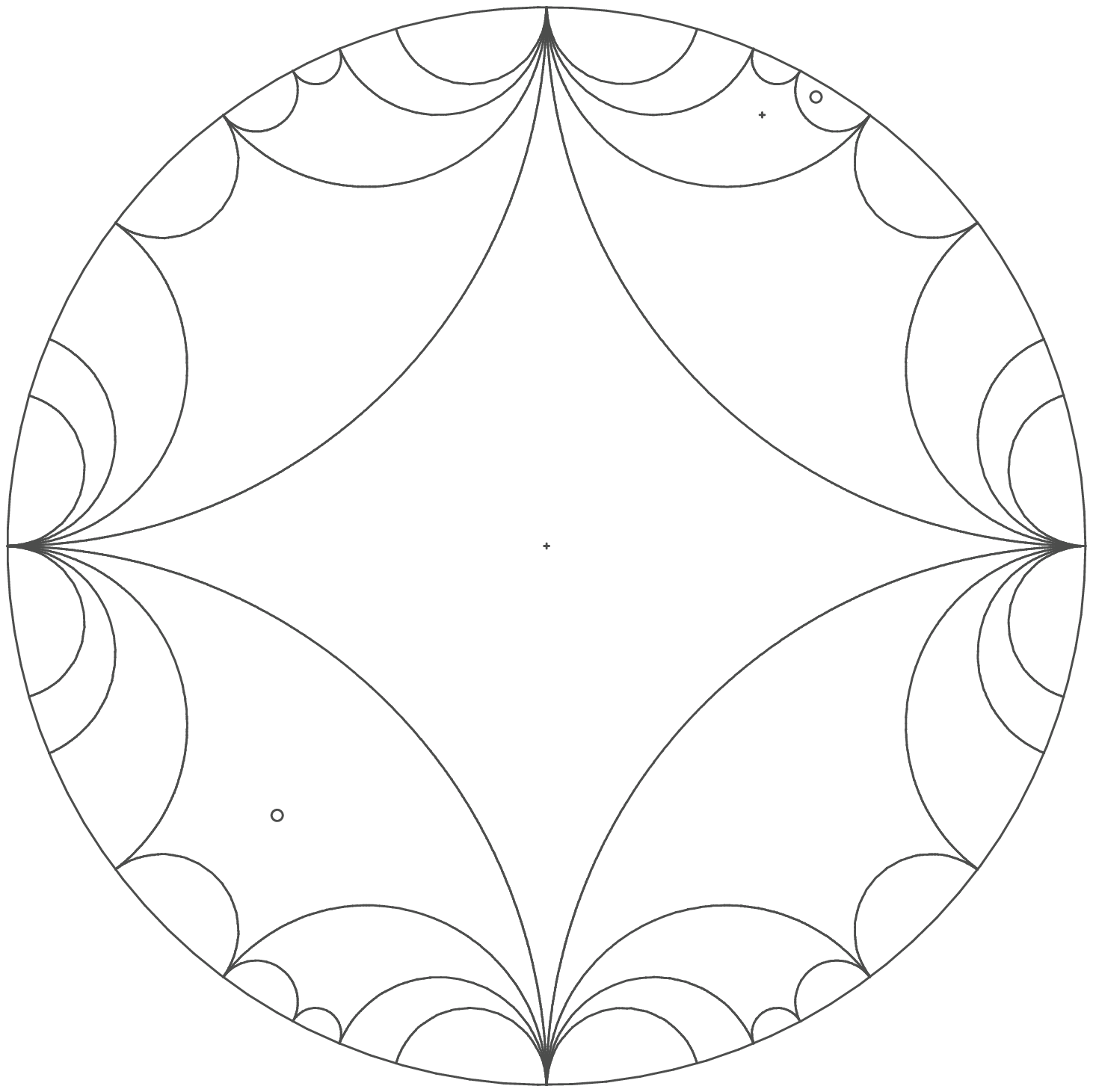}\\
Figure 2
\end{center}

\subsection{Non-isotopic diffeomorphisms}
As described in the introduction, there are many examples of pairs of manifolds
$X_0,X_1$ with $X_0 \# \cp^2 \cong X_1\# \cp^2$, where the diffeomorphism may
be assumed to carry the homology class of the $(+1)$-sphere $S_0$ to that of
$S_1$.  Adding two copies of $\cpbar$ gives $Z$, supporting two homotopic
diffeomorphisms $f_0, f_1$, which are precisely of the type discussed above. 
Let $f = f_1 \circ f_0^{-1}$.
\begin{theorem}
$D_Z(f)$ is a well-defined integer invariant and is given by
$$
D_Z(f) = D_Z(f_1) - D_Z(f_0) = -4(D_{X_1} -D_{X_0}).
$$
In particular, if $D_{X_0} \neq D_{X_1}$, then $f$ is a diffeomorphism which is
topologically isotopic to the identity, but not smoothly so.
\end{theorem}
\begin{proof}
Because $S_0$ and $S_1$ are homologous, $f_0$ and $f_1$ induce the same map on
homology, and hence $f$ induces the identity map on homology.  It
follows~\cite{cochran-habegger:homotopy,quinn:isotopy} that $f$ is homotopic and
in fact topologically isotopic to the identity of $Z$.  Since $f$ induces the
identity on $H_2(Z)$, it follows that $D_Z(f) \in \Z$ is defined. 
Lemma~\ref{additive}, together with the computation in Lemma~\ref{wallcross}, gives
the calculation of $D_Z(f)$.  If $D_{X_0}$ and $D_{X_1}$ are different, then it
follows from Lemma~\ref{isotopy} that $f$ is not isotopic to the identity.
\end{proof}

\subsection{The Seiberg-Witten version}\label{swversion}
It seems worth asking whether invariants based on the Seiberg-Witten equations
can be used to detect the non-triviality of a diffeomorphism.   Because the
diffeomorphism $f^{\Sigma^+}\circ f^{\Sigma^-}$ acts non-trivially on the set
of $\spinc$ structures, even up to the action of $J$, the current methods do
not give rise to a useful invariant.  However, using the ideas discussed in
section~\ref{swsec} we can in some circumstances define a mod $2$ invariant of
a diffeomorphism built from a single reflection, for example $f^{\Sigma^+}$. 
(Of course the hope would be to get different answers by taking $X$ to
be the manifolds $X_0$ and $X_1$.)  In order to do this, we need to assume that
the manifold $X$ has a $\spinc$-structure
$\ptilde_X\to X$ with $c_1(\ptilde_X) = 0$ and $\dim(\msw(\ptilde_X)) = 0$.
This hypothesis implies that $X$ is a spin manifold; using Rohlin's theorem,
and the formula for the dimension of $\msw(\ptilde_X)$, one readily deduces
that $b_+^2(X) \equiv -1 \pmod 4$.

Let $\ptilde_Z$ be the $\spinc$-structure on $Z$ which is a connected sum
of $\ptilde_X$ with the $\spinc$-structure on $N$ with $c_1$ Poincar\'e dual to
$\Sigma^+$. Then  $\ptilde_Z$ pulls back (via $f_0^{\Sigma^+}$) to
$J(\ptilde_Z)$, and so by Theorem~\ref{swbar} there is an invariant
$\SW(f,\ptilde_Z)$.  Since $b_+^2(Z) \equiv 0 \pmod 4$, the quantity $\epsilon$
appearing in that theorem is $1$, and so $\SW(f,\ptilde_Z)$ is only defined
modulo $2$.  In any event, we can still calculate $\SW(f,\ptilde_Z)$ by the
wall-crossing technique.  The details of the Kuranishi picture as the wall is
crossed are slightly different and may be deduced from the calculation
in~\cite{kronheimer-mrowka:thom}.   The result (whose proof we omit) should be
compared with Theorem~\ref{singlecalc}.
\begin{theorem}\label{swcalc}
Let $\ptilde_X\to X$ be a $\spinc$-structure with $c_1(\ptilde_X) = 0$ and
$\dim(\msw(\ptilde_X)) = 0$. Then the invariant $\SW(f,\ptilde_Z)$ described
above is equal (modulo $2$) to $\SW(\ptilde_X)$.
\end{theorem}

The rather unfortunate punchline to this story is that the mod $2$ invariant
$\SW(f,\ptilde_Z)$ does not depend on $X$ in any interesting way.  This follows
from a theorem of Morgan and Szab\'o~\cite{morgan-szabo:mod2} which states that
the invariant $\SW(\ptilde_X)$ for the sort of $X$ we considered is determined
(mod $2$) simply by $b_+^2$.  In fact, if $b_+^2(X) = 3$, then the
Seiberg-Witten invariant is odd, and if $b_+^2(X) > 3$, then the
Seiberg-Witten invariant is even.  So as it stands, we do not know how to use a
Seiberg-Witten invariant to detect non-isotopic diffeomorphisms.

\subsection{Discussion}
The construction of the examples seems in some respects more complicated than
necessary.   Before the addition of the two copies of $\cpbar$, we had the
manifold $Z'$, containing two spheres $S_0,S_1$ of self-intersection $+1$. 
Associated to each of these is a diffeomorphism $f_i'$, where the induced maps
in homology now satisfy $(f_i)_*(x) = x - 2(x\cdot S_i)S_i$.  It is natural to
conjecture that these diffeomorphisms are not isotopic  (As before, they could
certainly be arranged to be homotopic) and indeed that was our original
conception.  However, the construction given in this paper does not seem to
provide invariants of the $f_i'$ which would encode the difference in Donaldson
(or Seiberg-Witten) invariants of the manifolds $X_i$.

One aspect of this is related to the other complication in the examples, namely
that the diffeomorphisms $f_0$ and $f_1$ are not reflections, but are instead
compositions of such.  The calculation of the action on orientation of the
moduli space shows that a single reflection reverses orientations, and so leads
only to a $\Z_2$ invariant, which we have seen does not contain much useful
information about isotopy of diffeomorphisms.  This calculational fact has a
more conceptual basis, which can be seen if one interprets our invariants as
giving rise to cocycles of some sort in
$H^1(\met^*(Z)/\diff^+(Z))$.   A diffeomorphism of finite order preserves a
metric on $Z$, and hence the action on $\met(Z)$ has fixed points.  One gets a
better quotient space considering
$\met^*(Z)$, which consists of metrics with trivial isometry group.   One would
expect that if $f$ had order $n$, then the count of solutions around a loop in
$\met^*(Z)$ linking the fixed point would be well-defined only modulo
$n$.   If one views the examples in this light, then it is clear that one
should only have a mod $2$ invariant associated to a single reflection.  The
composition of two reflections (as used in the examples) has infinite order,
circumventing this difficulty.

An interpretation of our invariants as cohomology classes is the proper
setting for the generalization to invariants defined by families of dimension
greater than one.  We plan to return to this in a future paper.

\end{document}